\documentclass{article}
\usepackage{graphicx}
\usepackage{amsmath}

\newtheorem{thm}{Theorem}[section]

\newtheorem{prop}[thm]{Proposition}

\newtheorem{ex}{Example}[section]

\newcounter{alphthm}
\newcommand{\be}{\be}
\newcommand{\ee}{\end{equation}}
\newcommand{\ben}{\begin{enumerate}}
\newcommand{\een}{\end{enumerate}}
\newcommand{\pa}{{\partial}}

\newcommand{\pxi}{ {\pa \over \pa x^i}}

\def\beq{\begin{equation}} 
\def\eeq{\end{equation}}

\newcommand{\qed}{\hspace*{\fill}Q.E.D.}  

\title{Time-Optimal solutions of Parallel Navigation and Finsler geodesics}
\author{M. Rafie-Rad\footnote{Email address:\emph{m.rafiei.rad@gmail.com}}
\\
\small{Department of Mathematics, Faculty of Sciences,}\\
\small{Mazandaran University, Bablosar, Iran.}}

\begin{document}

\maketitle

\begin{abstract}
A geometric approach to kinematics in control theory is illustrated. A non-linear control system is derived for the problem and the Pontryagin maximum principle is used to find the time-optimal trajectories of the Parallel navigation. The time-optimal trajectories of the Parallel navigation are characterized through a geometric formulation. It is notable that the approach has the advantages using feedback.
\end{abstract}
\section{Introduction}
The historical development of what became the Calculus of Variations is closely linked to certain
minimization principles in the majority subjects in mechanics, namely, {\it the principle of least distance}, {\it the principle of least time} and ultimately, {\it the principle of least action} \cite{Lanc}. To understand solution of the well-known {\it brachistochrone} problem, (i.e finding a curve from point $A$ to point $B$ along which a free-sliding particle will descend more quickly than on any other $AB$-curve), we are led through {\it Fermat's principle of least time}: light always takes a path that minimizes travel time.\\
The {\it Parallel navigation}, or briefly P-navigation, is a quiet old problem and has been studied using several techniques from the viewpoints of kinematics and dynamics in optimal control theory \cite{Shney}. The application of Finsler geometry in Physics, seismology  and Biology is a subject of numerous papers such as \cite{ACH}, \cite{AnBoSl},\cite{AIM}, \cite{BR}, \cite{Mat}, \cite{Randers}, \cite{ShiSab}, \cite{YN}, etc. Let $O$ be the origin of an inertial reference frame of coordinates (FOC). The positions of $M$ and $T$ in this (FOC) are given by the vectors ${\bf r}_M = OM$ and ${\bf r}_T = OT$, respectively. In two-point guidance systems, the vector ${\bf r}={\bf r}_T -{\bf r}_M$ is conventionally called the {\it range}. Its time derivative $\dot{\bf r}=\dot{\bf r}_T -\dot{\bf r}_M={\bf v}_T-{\bf v}_M$ is {\it the relative velocity} between the two objects, and ${\bf v}_T$ and ${\bf v}_M$ are the velocities of $T$ and $M$, respectively. W always denote the vectors by bold face and their norms will be shown by the same normal letter. As an application, it is notable for mariners wishing to rendez-vous each other at sea. $M$ could be a boat and $T$, a tanker with fuel for it (or vice-versa). Or, back in history, $T$ could be a merchantman and $M$ a pirate ship. This rule assumes, of course, constant speeds. Thus, in most realistic cases, $v_T$ and $v_M$ are supposed to be constant. However, it is easy to extend the theory if they are not constant.  The {\it closing velocity}, a term often used in the study of
guidance, is simply ${\bf v}_C=-\dot{\bf r}$. Notice that, we wish to study the kinematics of P-navigation in a relative (FOC) rather than a absolute one, i.e., we shall seek the location of $M$ in a (FOC) attached to $T$. Thus, a trajectory
in the relative (FOC) shows the situation as seen by an observer located at $T$. As the special cases, we assume that $M=R^3$ or $M=R^2$. In reality, the velocity ${\bf v}_T$ and ${\bf r}_T$ can be detected and reported at any ${\bf r}$ by a grounded radar. Suppose that $\delta({\bf r})$ be the angle between ${\bf v}_M$ and $MT$ and given any $\delta$, there is Finsler metric $F$ given by:
\begin{equation}
\label{Mats}
F({\bf r},{\bf v},\delta)=\frac{|{\bf v}|^2}{v_M\cos\delta|{\bf v}|-\langle{\bf v},{\bf v}_T\rangle},
\end{equation}
where, $|.|$ denotes the Riemannian norm on $M$. A {\it solution} of the described P-navigation is a curve $({\bf r}(t),\delta(t))$ such that respects the required constraints on velocities.
\begin{thm}
\label{mainthm1}
Given any solution $({\bf r},\delta)$ of parallel navigation, the curve ${\bf r}$ can be reparametrized so that it satisfies $F({\bf r}(t),{\bf v}(t),\delta(t))=1$.
\end{thm}
The \textit{indicatrix} $S({\bf r},\delta)$ of the metric (\ref{Mats}) is the set of unit tangent vectors ${\bf v}$ with respect to (\ref{Mats}) which is defined by $S({\bf r},\delta)=\{{\bf v}\in T_{\bf r}M\ |\ F({\bf r},{\bf v},\delta)=1\}$. Following Theorem \ref{mainthm1}, at any time $t$ we have $\dot{\bf r}={\bf v}\in S({\bf r},\delta)$. Hence, at any time $t$, there is a unit vector $f({\bf r},\delta)\in S({\bf r},\delta)$ such that $\dot{\bf r}={\bf v}=f({\bf r},\delta)$.\\
Control problems typically concern finding a (not necessarily unique) control law
$\delta(.)$ , which transfers the system in finite time from a given initial state $x_i={\bf r}(0)$ , to a given
final state $x_f={\bf r}(t_f)$. This transition is to occur along an admissible path, i.e. ${\bf r}(.)$ and respects all kinematic constraints imposed on it. Let us consider it as
\begin{equation}
\label{Control system}
\dot{\bf r}=f({\bf r},\delta).
\end{equation}
We further assume that $\delta(.)$ is admissible, i.e. is piecewise continuous and belongs to
${\cal U}$ , the admissible control space. Let there now be a rule which assigns a unique, real-valued number to each of these transfers. Such a rule can be viewed as the transition cost between $x_i$ and $x_f$ along
an admissible path, completely specified by $\delta(.)$. The Optimal control concerns
specifying this rule and thereby providing a systematic method for selecting the "best",
or "optimal" control law, according to some prescribed cost functional. One can find an analogue discussion in \cite{BR}, to calculate the travel-time along the trajectories of the so called \textit{Pure pursuit navigation}. Here, the P-navigation optimal control problem can be founded by the cost function $C({\bf r},\delta)=F({\bf r},\dot{\bf r},\delta)$ and has the following form
\begin{equation}
\label{min}
\textrm{minimize}\int_0^{t_f}C({\bf r},\delta)dt,
\end{equation}
where, $t_f\in(0,\infty)$ is the final time which is going to be optimized. From everyday experience we know that collision courses need not be straight lines if $T$ changes its speed or direction; so what is exactly the collision course? It may be curved in some sense. One of our goal in this paper is to make known the best collision course.
\begin{thm}
\label{mainthm2}
Given any time-optimal solution $({\bf r},\delta)$ of P-navigation, the curve ${\bf r}$ is a geodesic of the Finsler metric (\ref{Mats}).
\end{thm}
The trajectory ${\bf r}_M$ can be obtained ${\bf r}_M={\bf r}_T-{\bf r}$ when ${\bf r}$ is known. One can freely consider ${\bf v}_M$ and ${\bf v}_T$ as vector fields alon ${\bf r}$. Now, let $\frac{\nabla }{\ dt}$ be the covariant derivative defined for any vector field $Y$ along ${\bf r}$ defined by
\[
\frac{\nabla Y^i}{dt}:=\frac{dY^i}{dt}+G^i_{jk}({\bf r},\dot{\bf r},\delta)Y^jY^k,
\]
where, $G^i_{jk}$ are the connection coefficients of Berwald connection associated to the Finsler metric (\ref{Mats}).
As a result of Theorem \ref{mainthm2}, we can mention the following result:
\begin{thm}
\label{mainthm3}
The time-optimal trajectory ${\bf r}_M$ of P-navigation satisfies the following second order ODE:
\[
\ddot{\bf r}_M^i+G^i_{jk}({\bf r},{\bf v},\delta){\bf v}_M^j{\bf v}_M^k=\frac{\nabla{\bf v}_T^i}{dt},\ \ \ \ i=1,...,n.
\]
\end{thm}
Our approach is closely related with subjects such as non-holonomic mechanics, sub-Finslerian geometries, see for a deeper sight \cite{LM} and \cite{Bl}. One may find various techniques in missile guidance and control in \cite{Shney}.
\section{Preliminaries}
Let $M$ be a n-dimensional $ C^\infty$ manifold.  $T_x M $ denotes  the tangent space of M at $x$. The tangent bundle of $M$ is the union of tangent spaces $TM:=\cup _{x \in M} T_x M$.  We will denote  the elements of $TM$ by $(x, y)$ where $y\in T_xM$. Let $TM_0 = TM\setminus \{ 0 \}.$ The natural projection $\pi: TM_0 \rightarrow M$ is given by $\pi (x,y):= x$.\\
 A {\it Finsler metric} on  $M$ is a function $ F:TM \rightarrow [0,\infty )$ with the following properties; (i) $F$ is
$C^\infty$ on $TM_0$,\  (ii) $F$ is positively 1-homogeneous on the fibers of tangent bundle $TM$,  and  (iii) the $y$-Hessian of $\frac{1}{2}F^{2}$
with elements $ g_{ij}(x,y):=\frac{1}{2}[F^2(x,y)]_{y^iy^j} $ is positive definite on $TM_0$.  The pair  $(M,F)$ is then called a {\it Finsler space}. The Riemannian metrics are special Finsler metrics. Traditionally, a Riemannian metric is denoted by $ a_{ij}(x) dx^i \otimes dx^j$. It is a family of inner products on tangent spaces. Let
 $\alpha(x,y):= \sqrt{g_{ij}(x)y^iy^j}$,
${\bf y}= y^i\pxi|_x\in T_xM$. $\alpha$ is  a family of Euclidean norms on tangent spaces. Throughout this
paper, we also denote   a Riemannian metric by  $\alpha=\sqrt{a_{ij}(x)y^iy^j}$.

An $(\alpha, \beta)$-metric is a scalar function on $TM$ defined by $F:=\Phi(\frac{\beta}{\alpha})\alpha$,
where  $\phi=\phi(s)$ is a $C^\infty$ on $(-b_0, b_0)$ with certain regularity. $\alpha=\sqrt{a_{ij}(x)y^iy^j}$ is a Riemannian metric and $\beta =b_i(x)y^i$ is a 1-form on a manifold $M$. One may find another important class of $(\alpha,\beta)$-metrics in \cite{Shib}. The {\it Randers} and {\it Matsumoto} metrics are special  $(\alpha, \beta)$-metrics defined by $\Phi=1+s$ and $\Phi=\frac{1}{1-s}$, respectively, i.e, $ F = \alpha +\beta$ and $F=\frac{\alpha^2}{\alpha-\beta}$. Randers metrics were introduced by Randers in 1941 \cite{Randers} in the context of general relativity. In \cite{BaRoSh}, applying Fermat's principle, the authors proved that the time-optimal solutions of the well-known Zermelo's navigation-moving that is the motion of a vehicle equipped with an engine with a fixed power
output in presence of a wind current-are actually the geodesics of a Randers metric. M. Matsumoto gave an exact formulation of a Finsler surface to measuring the time on the slope of a hill and introduced the Matsumoto metrics in \cite{Mat}, see also \cite{ShiSab}.

A Lagrangian on the manifold $M$ is a mapping $L:TM\longrightarrow R$ which is smooth on $TM_0$. A Lagrangian is said to be {\it regular} if it has non-degenerate $y$-Hessian on $TM_0$. Thus, given a Finsler metric $F$, the function $L=\frac{F^2}{2}$ is a regular Lagrangian. A large area of applicability of this geometry is
suggested by the connections to Biology, Mechanics, and Physics and also by its
general setting as a generalization of Finsler and Riemannian geometries \cite{MiAnas}. For every smooth curve $c:[a,b]\longrightarrow R$, the extremal curves of the action integral given by
\begin{equation}
\label{action}
I(c)=\int_a^bL(c(t),\dot{c}(t))dt,
\end{equation}
are characterized locally by the {\it Euler-Lagrange equations} given as follows:
\begin{equation}
\label{Euler-Lagrange}
\frac{d}{dt}\frac{\partial L}{\partial\dot{x}^i}-\frac{\partial L}{\partial x^i}=0,
\end{equation}
where, $x^i(t)$ is a local coordinate expression of $c$. The extremal curves of the action integral (\ref{action}) are usually called {\it the geodesics of L}. In \cite{ACH} it is shown that the Lagrangian and Finslerian approaches are projectively the same.\\
Given a Finsler manifold $(M,F)$, a globally defined vector field $G$ is
induced by $F$ on $TM_0$, which in a standard coordinate $(x^i,y^i)$
for $TM_0$ is given by
$
G=y^i {{\partial} \over {\partial x^i}}-2G^i(x,y){{\partial} \over
{\partial y^i}},
$
where $G^i(x,y)$ are local functions on $TM_0$ satisfying $G^i(x,\lambda y)=\lambda^2 G^i(x,y)\,\,\, ,\lambda>0$, see \cite{Sh}. G is called the  associated {\it {spray}} to $(M,F)$. In
local coordinates, a curve $c(t)$ is a geodesic of $F$ if and only if its coordinates $(c^i(t))$ satisfy $ \ddot c^i+2G^i(c,\dot c)=0$.

\subsection{The kinematics of Parallel navigation}
\begin{figure}
  \hspace{3cm}
  \includegraphics[width=7cm]{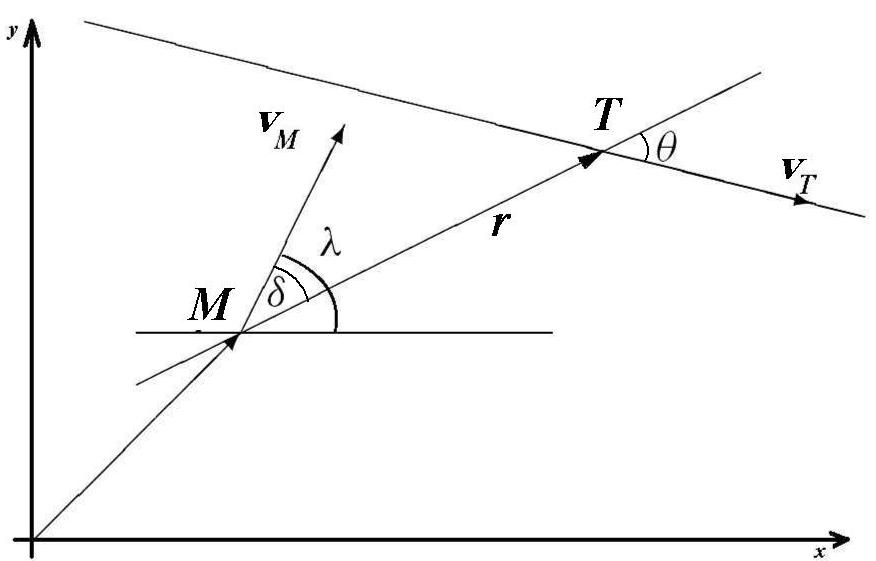},
  \caption{\textit{The range ${\bf r}$, the velocity vectors ${\bf v}_M$ and ${\bf v}_T$.}}\label{F1}
\end{figure}

We shall refer to the target as $T$ and to the pursuer as $M$ and their velocities as $v_M$ and $v_T$, respectively. To begin, we set up a coordinate system called reference frame of coordinates,
in which the pursuer is initially located at the origin $O$. When considering planar motion we shall use Cartesian
coordinates $(x, y)$ or $(x, z)$, and the angles will be positive if measured counterclockwise. The ray that starts at the pursuer $M$ and is directed at the target $T$ along the positive sense of ${\bf r}$ is called \textit{the line of sight} (LOS). The \textit{parallel navigation} geometrical rule,has been known since antiquity, mostly by mariners. According to this rule, the direction of the line of sight, $MT$, is kept constant relative to inertial space, i.e., the LOS is kept parallel to the initial LOS. In three-dimensional vector terminology, the rule is very concisely stated as ${\bf r}\times\dot{\bf r}=0$. Suppose that $\theta$ and $\lambda$ denote, respectively, the angles between ${\bf v}_T$ and ${\bf v}_M$ and, ${\bf v}_M$ and the
horizontal axis (Figure \ref{F1}).

Let us put $r=|{\bf r}|$. The basic rule for moving of the pursuer is presented by the following two equations \cite{Shney}:
\begin{eqnarray}
\label{PP law 1}
\dot{r}&=&v_T\cos\theta-v_M\cos\delta,\\
r\dot{\lambda}&=&v_T\sin\theta-v_M\sin\delta.\label{PP law 2}
\end{eqnarray}
Notice that, in a planar framework, ${\bf v}_M$ , ${\bf v}_T$ and ${\bf r}$ being on the same (fixed) plane by definition, therefore, the parallel navigation geometrical rule can be restated as $\dot{\lambda}=0$. The requirement $\langle{\bf r},{\bf v}\rangle<0$ must be added in order to ensure that $M$ should approach $T$ not recede from it. In this case, we have $\dot{r}<0$, that is $v_T\cos\theta<v_M\cos\delta$. Let us denote the projection of any vector ${\bf v}_T$ on ${\bf v}$ by $Proj_{{\bf v}}{\bf v}_T$. A \textit{solution} of the described P-navigation is a curve $({\bf r}(t),\delta(t))$ such that respects the equations (\ref{PP law 1}) and (\ref{PP law 2}). By \textit{the trajectory} of P-navigation, we mean a curve ${\bf r}(t)$ such that $({\bf r}(t),\delta(t))$ is a solution, for some control $\delta$.

\begin{figure}
  \hspace{3cm}
  \includegraphics[width=5cm]{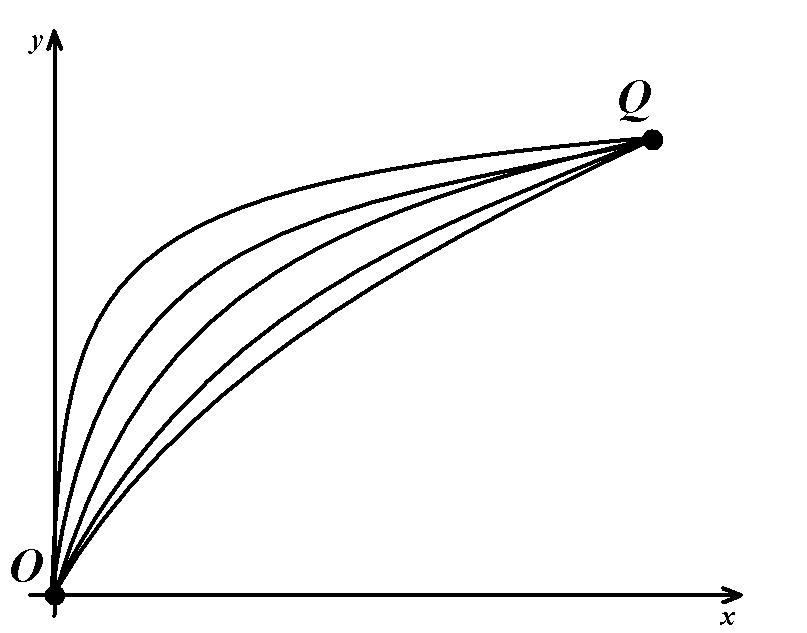},
  \caption{\textit{Some possible ranges initiated at the point $Q$.}}\label{F2}
\end{figure}

\begin{figure}
  \hspace{2cm}
  \includegraphics[width=8cm]{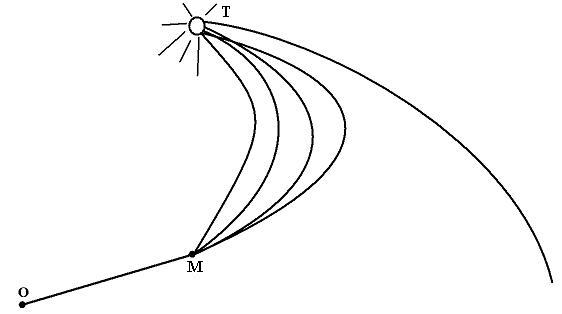},
  \caption{\textit{Schematic of exemplary collision courses for $M$.}}\label{F3}
\end{figure}

Initiating the process, we have ${\bf r}(0)={\bf r}_0$ which shows that, $M$ stands at a point with distance $r_0$ from $T$. Through the performance, $r$ decreases by time and hence, $M$ approaches $T$. Therefore, ${\bf r}$ tends to the origin $O$ and $M$ hits $T$ when ${\bf r}(t_f)=0$, (Figure \ref{F2}). It follows that, P-navigation trajectories are characterized by a curve ${\bf r}$ joining $Q={\bf r}_0$ to the origin $O$ (Figure \ref{F3}). It is of our interests to find the best $QO$-trajectory. More precisely, the problem is \textit{to find a curve from point $Q$ to point $O$ along which a particle will descend more quickly than on any other $QO$-curve of P-navigation.} In this way, the problem somehow resembles to a brachistochrone problem.
\bigskip

\section{The optimal control theory.}
A \textit{control system} of ordinary differential equations is a family of differential equations in normal form $\frac{d{\bf r}^i}{dt}= f^i({\bf r}, \delta)$, where ${\bf r}^i$ are called \textit{state variables}, $t$ is the \textit{parameter of evolution} (usually the time) and $\delta^a$ are the \textit{controls}. Geometrically, it can be regarded as a fibred mapping $X:U\longrightarrow TM$, from a control fiber
bundle $(U,\eta,M)$ over the state manifold $M$ to the tangent bundle $(TM,\pi,M)$, see \cite{NvS}. Using
local coordinates $({\bf r}^i),\ i=1,...,n$ in $M$, adapted coordinates $({\bf r}^i, \delta^a),\ a=1,...,k$ in $U$, and natural coordinates $({\bf r}^i, {\bf v}^i)$ in $TM$, the coordinate expression for $X$ is $X({\bf r}, \delta) = f^i({\bf r},\delta )\frac{\partial}{\partial {\bf r}^i}$ ,
or ${\bf v}^i = f^i({\bf r}, \delta)$, the family of control equations. \textit{Admissible} curves of the control system are curves $\gamma: I\subset R\longrightarrow U$ such that $(\eta o\gamma)^c = X o\gamma$, where $^c$ denotes the natural lifting to $TM$ of a curve  in $M$. Interested readers are advised to see \cite{NvS} for getting familiar to the geometry of control systems. In Optimal Control Theory, a {\it cost functional} ${\cal C}(\gamma) = \int C({\bf r}(t), \delta(t))dt$ is given and the goal is to obtain admissible curves of the control system, satisfying some boundary conditions (e.g. $x_i={\bf r}(0)$, $x_f={\bf r}(t_f)$) and minimizing the cost functional. It is therefore a Classical Variational problem with non-integrable constraints defined
by the control equations. Pontryagin maximum principle \cite{Po} provides a set of necessary conditions for a solution $({\bf r}(t),\hat{\delta}(t))$ to be optimal; introducing a Hamiltonian function
\begin{eqnarray*}
H({\bf r},{\bf p} , \delta) &:=&\langle {\bf p},X\rangle-C({\bf r},\delta)=  {\bf p}_if^i({\bf r}, \delta)-C({\bf r},\delta),\\
\hat{H}({\bf r},{\bf p})&:=&\underset{\delta}{\max}\ H({\bf r},{\bf p} , \delta).
\end{eqnarray*}
where the variables $({\bf p}_i)$ are momenta coordinates, the optimal curves $({\bf r}(t), \hat{\delta}(t))$ must satisfy the control system equations
\[
{\bf v}^i=\frac{\partial \hat{H}}{\partial {\bf p}^i}=f^i({\bf r}(t), \hat{\delta}(t))
\]
and there must exist a solution curve for the adjoint differential equations
\[
\frac{d{\bf p}_i}{dt}=-\frac{\partial \hat{H}}{\partial {\bf r}^i},
\]
Define the Lagrangian $L$ by $L({\bf r},{\bf v})={\bf p}_i{\bf v}^i-\hat{H}$.
Observe that we have the following relations
\[
\frac{d{\bf r}}{dt}= \frac{\partial \hat{H}}{\partial {\bf p}}={\bf v},\ \ \ \ \ \ \frac{d{\bf p}}{dt}=-\frac{\partial \hat{H}}{\partial {\bf r}}=\frac{\partial L}{\partial {\bf r}}, \ \ \ \ \ \frac{\partial \hat{H}}{\partial {\bf v}}={\bf p}-\frac{\partial L}{\partial {\bf v}}=0.
\]
From the above equations, it results the well-known Euler-Lagrange for $L$
\[
\frac{d}{dt}\frac{\partial L}{\partial {\bf v}}-\frac{\partial L}{\partial {\bf r}}=0.
\]
\begin{prop}
\label{Pontryagin}
\cite{Po} In order for $({\bf r}(t),\hat{\delta}(t))$ to be an optimal solution of (\ref{min}), the following are necessary conditions:\\
(a)  There exists a solution curve for the adjoint differential equations $$\frac{d{\bf p}_i}{dt}=-\frac{\partial \hat{H}}{\partial {\bf r}^i}.$$\\
(b)  $\hat{\delta}=\arg \ \underset{\delta}{\max} \ H ({\bf r},{\bf p},\delta),\ \ \ \ \forall t\in[0,t_f]$.\\
(c) $\hat{H}({\bf r},{\bf p}) = 0, \ \ \ \ \forall t\in[0,t_f]$.
\end{prop}

\section{Proof of Theorems.}
\subsection{Proof of Theorem \ref{mainthm1}}
Let $({\bf r}(t),\delta(t))$ be a pair of the curve ${\bf r}$ and a function $\delta(t)$. We are going to show that, if $({\bf r}(t),\delta(t))$ be a solution of P-navigation, then ${\bf t}(t)$ must be reparametrized so that we we have $F({\bf r}(t), \dot{\bf r}(t),\delta(t))=1$. We notice that, in P-navigation, ${\bf r}$ and ${\bf v}$ are collinear and $\dot{r}<0$, hence we have
\[
\dot{r}=\frac{\langle{\bf r},{\bf v}\rangle}{r}=\pm|Proj_{\bf r} {\bf v}|=\pm|Proj_{\bf v} {\bf v}|=-|{\bf v}|.
\]
Now, we summarize (\ref{PP law 1}) in the following relation
\[
|{\bf v}|=v_M\cos\delta-\frac{\langle{\bf v}_T,{\bf v}\rangle}{|{\bf v}|}.\\
\]
After simplification, we obtain the following equation
\[
F({\bf r},{\bf v},\delta)=\frac{|{\bf v}|^2}{v_M\cos\delta|{\bf v}|-\langle{\bf v}_T,{\bf v}\rangle}=1.
\]
\qed
\subsection{Proof of Theorem \ref{mainthm2}}
Following Theorem \ref{mainthm1}, at any time $t$ we have $\dot{\bf r}={\bf v}\in S({\bf r},\delta)$. Hence, at any time $t$, there is a unit vector $X({\bf r},\delta)\in S({\bf r},\delta)$ such that $\dot{\bf r}={\bf v}=X({\bf r},\delta)$.
Consider the unit canonical vector field $\ell({\bf r},\dot{\bf r},\delta)=\frac{\dot{\bf r}}{F({\bf r},\dot{\bf r},\delta)}$. We notice that, in P-navigation framework, we always assume that ${\bf r}$ and $\dot{\bf r}$ are collinear and hence, one can understand $\ell$ as a function of ${\bf r}$ and $\delta$, as well. It follows that, given any trajectory ${\bf r}$ of P-navigation, $X$ is given by $X({\bf r},\delta)=\ell({\bf r}, \dot{\bf r},\delta)$. Therefore, it is clear that,
\begin{eqnarray*}
\langle {\bf p},X\rangle &=&p_if^i({\bf r},\delta)=p_i\ell^i({\bf r},\dot{\bf r},\delta)=F({\bf r},\dot{\bf r},\delta),\\
\langle {\bf p},{\bf v}\rangle &=&p_i{\bf v}^i=F^2({\bf r},\dot{\bf r},\delta).
\end{eqnarray*}
Now, we return to the control system of P-navigation given by (\ref{Control system}) with the cost functional $C({\bf r},\delta)=F({\bf r},\dot{\bf r},\delta)$.
It is easy to verify that, $H=0$, $\hat{H}=0$ and one may consider $\hat{\delta}$ as any possible control law. The conditions of Proposition \ref{Pontryagin} holds as well and the Lagrangian $L_{\hat{\delta}}=\langle{\bf p},{\bf v}\rangle-\hat{H}$ is obtained as
\[
L_{\hat{\delta}}({\bf r},\dot{\bf r})=F^2({\bf r},\dot{\bf r},\hat{\delta}).
\]
Therefore, based on Pontryagin maximum principle, the optimal trajectories ${\bf r}(t)$ are geodesics of the Lagrangian $L_{\hat{\delta}}$. Clearly, they are geodesics of the Finsler metric $F({\bf r},\dot{\bf r},\delta)$.

Now, consider the control-parametric family of Finsler metrics defined by $F_\delta({\bf r},\dot{\bf r}):=F({\bf r},\dot{\bf r},\delta)$. Let ${\cal L}_\delta(\gamma)=\int_0^{t_f}F_\delta(\gamma,\dot{\gamma})dt$ be the length of any admissible curve $\gamma(t)$ on $(M,F_\delta)$. A simple calculation gives the following inequality:
\[
F_0({\bf r},\dot{\bf r})\leq F_\delta({\bf r},\dot{\bf r}),\ \ \ \textrm{for all possible controls}\ \delta.
\]
From that, it follows that the functional ${\cal L}_\delta(\gamma)$ takes its minimum at  $\delta=0$, that is
\[
{\cal L}_0(\gamma)\leq{\cal L}_\delta(\gamma),\ \ \ \textrm{for all possible controls}\ \delta.
\]
Therefore, to find a time-optimal solution, one should minimize the cost functional ${\cal C}(\gamma)=\int F_0(\gamma,\dot{\gamma})dt$ and this leads us to obtain it as a geodesic of $F_0$.
\qed
\begin{thm}
The time-optimal trajectory of P-navigation is a geodesic ${\bf r}(t)$ of the Finsler metric $F_0=\frac{|{\bf v}|^2}{v_M|{\bf v}|-\langle{\bf v}_T,{\bf v}\rangle}$.
\end{thm}
However, given any control law, one may obtain a geodesic of the metric $F_\delta$ as the time-optimal trajectory. As a remark, we quote that the target $T$ may not be reachable by the control $\delta=0$.
\begin{ex}
(Case of plane nonmaneuvering target.)\textrm{\textsf{\textrm{ The target $T$ is said to be \textit{nonmaneuvering} if ${\bf a}_T=0$. In this case, $T$ moves on a straghit line at velocity $v_T$ in the direction with a constant angle $\theta_0$ if measured counterclockwise, see Figure \ref{F4}. Let us suppose ${\bf v}_T(x^1,x^2)=v_T\{\cos\theta_0\frac{\partial}{\partial x^1}+\sin\theta_0\frac{\partial}{\partial x^2}\}$. Thus, from (\ref{PP law 2}), it follows that $\delta=\sin^{-1}(\frac{\sin\theta_0}{K})$, where, $K$ is the {\it velocity ratio} $K = \frac{v_M}{v_T}$. Then, $\delta$ is a constant say $\delta_0$. Moreover, ${\bf v}_T$ is a parallel vector field and then $F_\delta$ is a Minkowski metric and is flat. Thus, it geodesics are straight lines. We obtain ${\bf r}(t)={\bf r}_0+t{\bf v}_0$. But, from (\ref{PP law 1}), we have $|{\bf v}|=|{\bf v}_0|=v_M\cos\delta_0-v_T\cos\theta_0$. Intercept occur when we have ${\bf r}(t_f)=0$, thus, the total flight time $t_f$ is obtained by $$t_f=\frac{r_0}{v_M\cos\delta_0-v_T\cos\theta_0}=\frac{r_0}{v_T(K\cos\delta_0-\cos\theta_0)}$$ and the total range of $M$ equals $r_0$ which is the shortest curve joining ${\bf r}_0$ to the origin $O$. }}}
\begin{figure}
  \hspace{2cm}
  \includegraphics[width=7cm]{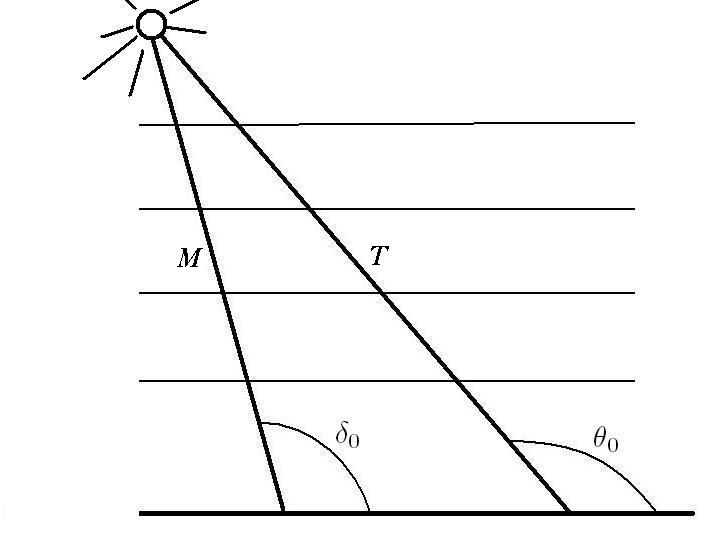},
  \caption{\textit{Collision course for a target moving on a straight line at a direction with a constant angle $\theta_0$.}}\label{F4}
\end{figure}
\end{ex}

\end{document}